\documentclass[11pt,a4paper,reqno]{amsart}
\usepackage{graphicx} 

\usepackage{amsmath}
\usepackage[foot]{amsaddr}

\usepackage{enumitem}
\usepackage{xcolor}
\usepackage{amssymb}
\usepackage{amsthm}
\usepackage{array}
\newcolumntype{M}[1]{>{\centering\arraybackslash}m{#1}}

\usepackage[backend=bibtex]{biblatex}
\usepackage{float}

\usepackage{cuted}
\usepackage{color}
\usepackage{graphicx}

\usepackage{tikz}
\usetikzlibrary{shapes.geometric}
\usetikzlibrary{arrows.meta,arrows}

\usepackage{tabularx}

\usepackage{comment}

\usepackage[hidelinks, colorlinks=false]{hyperref}

\usepackage{cleveref}

\usepackage[normalem]{ulem}

\newtheorem{thm}{Theorem}

\newtheorem{definition}[thm]{Definition}
\newtheorem{example}[thm]{Example}

\newcommand{\R}{\mathbb{R}}
\newcommand{\N}{\mathbb{N}}

\newcommand{\scp}[2]{ \left \langle #1, #2 \right\rangle }
\newcommand{\n}[1]{\left\Vert #1\right\Vert}

\newcommand{\eqb}{\begin{flalign}}
\newcommand{\eqe}{\end{flalign}}

\newcommand{\Lra}{\Leftrightarrow}

\newcommand{\platz}{\;\;\;\;\;}

\allowdisplaybreaks

\definecolor{mygrey}{RGB}{220, 220, 220}
\definecolor{concrete}{RGB}{136, 150, 136}
\definecolor{steel}{RGB}{200, 200, 200}
\definecolor{copper}{RGB}{220, 127, 100}
\definecolor{redd}{RGB}{200, 50, 50}
\definecolor{yellow}{RGB}{255, 254, 235}
\definecolor{green}{RGB}{213, 243, 231}
\definecolor{green2}{RGB}{169, 209, 142}
\definecolor{lightgrey}{RGB}{240, 240, 240}
\definecolor{lightred}{RGB}{255, 244, 239}
\definecolor{lightlilac}{RGB}{239, 239, 255}
\definecolor{wood}{RGB}{133, 94, 66}
\definecolor{BackgroundGrey}{RGB}{210, 210, 210}

\begin{document}

\title[Spatial decay of perturbations in a wave equation with optimal boundary control]{A simulation study on spatial exponential decay of perturbations in a two-dimensional wave equation with optimal boundary/line control}
\thanks{B.O.\ extends his gratitude to the Thüringer Graduiertenförderung for funding his scholarship.}
\author{Benedikt Oppeneiger$^{1,2}$}
\address{$^1$Faculty of Mathematics, Chemnitz University of Technology, Germany}
\address{$^2$Optimization-based Control Group, Institute of Mathematics, Technische Universität Ilmenau, Germany.\\ Mail: \textsc{benedikt-florian.oppeneiger@tu-ilmenau.de}}
\date{\today}

\begin{abstract}
Recent results have shown that domain-uniform stabilizability and detectability imply spatial exponential decay of perturbations in optimally controlled hyperbolic PDEs on one-dimensional domains.
This domain-uniform stabilizability and detectability can be achieved only if the control domain is distributed over the whole spatial domain such that the distance between neighboring control intervals is bounded from above.
Motivated by these insights we investigate whether analogous effects can be observed for optimal boundary control in higher dimensions.
For this purpose we conduct a simulation study 
on a two-dimensional wave equation on expanding square domains which is driven by a localized perturbation of the initial displacement. 
We compare two control geometries: a control acting only on the outer boundary and a control acting on a regular grid of interior line interfaces combined with the boundary.
The problem is discretized by conforming \(P_1\) finite elements in space and the implicit midpoint rule in time. 
Instead of assembling the full space-time KKT system, we use a condensed formulation and solve the reduced optimality system by preconditioned conjugate gradients.
We show that domain-uniform spatial decay of the perturbation can only be observed in the scenario with a regular grid of line controls. This is because finite propagation velocity requires a uniform bound on the distance that a wave can travel without reaching the control domain.
\end{abstract}
\maketitle                   

\section{Introduction}
\label{Sec: Introduction}
\noindent Robustness is an important concern in the analysis of optimal control problems governed by partial differential equations. In general, this refers to the ability of a system to tolerate perturbations or uncertainties without exhibiting an excessive change in its behaviour. Specifically in the context of optimal control small changes in the data should not lead to disproportionate changes in the optimal state, adjoint state, or control.
 In many real-world scenarios perturbations may be confined to a small region. 
 Typical examples include a local voltage fluctuation in a power grid, a sudden braking maneouver of a small group of cars on a highway, or a toxic substance entering a drinking water network via a small leak in a pipeline. 
 In all of these cases a global propagation of the local perturbation causing phenomenons like a blackout, traffic congestion or large-scale water pollution is highly undesirable.

\noindent The fundamental idea of local perturbations only having a local impact has been the subject of substantial mathematical research across several different areas.
In optimal control, a closely related phenomenon is the turnpike property (see~\cite{Damm2014, Breiten2020, Faulwasser2022}).
There, an initial condition may be interpreted as a local in time perturbation of the optimal steady state. As a consequence of turnpike the distance of the optimal trajectory to this equilibrium decays over time.
These ideas were extended in~\cite{Schaller2020} and~\cite{Schaller2022} where an exponential temporal decay of more general local in time perturbations was proven for optimal control problem governed by general evolution equations.

\noindent Another locality principle appears in graph-structured optimization. In \cite{Shin2022}, it is shown that sensitivities of primal-dual solutions of graph-structured nonlinear programs may decay exponentially with the graph distance from the perturbed data. 
Related ideas are used in diffusing-horizon model predictive control \cite{Shin2023}, where local information is propagated only over a limited neighbourhood in the underlying graph. 
Such graph-structured optimization problems may have their origin in discretized PDE-constrained optimal control problems, where the network structure corresponds to the geometry of the computational mesh.

\noindent Interesting results in this context have also appeared in numerical linear algebra. 
In \cite{Demko1984} it was shown that under suitable invertibility and conditioning assumptions, the entries of the inverse of a band matrix decay exponentially away from the diagonal. 
Since banded and sparse matrices naturally arise from local discretizations of differential operators, such results can be viewed as a discrete precursor of spatial decay phenomena for PDEs. 
In the more recent work~\cite{Benzi2015} these estimates were refined.

\noindent For linear-quadratic PDE-constrained optimal control, spatial decay of perturbations was first studied in \cite{Goettlich2025} for elliptic and parabolic problems. 
In this framework, local perturbations of the data lead to exponentially localized perturbations of the optimal state, adjoint state, and control, provided that the associated KKT solution operators are bounded uniformly with respect to the size of the spatial domain.
To achieve this domain-uniform boundedness a particular stabilizability/detectability condition was introduced.
This property requires the existence of a bounded state feedback such that the closed loop operators are coercive with constants which are again uniform in the domain size.

\noindent This mechanism was later significantly generalized to general evolution equations, including hyperbolic problems, in \cite{Oppeneiger2026}. In contrast to the coercivity-based condition used for elliptic and parabolic equations, the latter work employs a semigroup-based notion of domain-uniform stabilizability and detectability. This condition requires the existence of uniformly bounded feedback operators such that the closed-loop semigroups decay exponentially with constants independent of the spatial domain. 
For several hyperbolic equations, this property was characterized by a uniform spatial distribution of the control and observation regions. 
Intuitively, due to the finite propagation velocity the distance between neighbouring controlled parts of the spatial domain must bounded from above, and the controlled sets inside prescribed subdomains must not degenerate. The rigorous characterization in \cite{Oppeneiger2026}, however, is restricted to one-dimensional domains and bounded distributed input and output operators. Boundary control and lower-dimensional control sets, such as line controls in two space dimensions, are therefore not covered by the existing theory.

\noindent The present paper is motivated by the conjecture that analogous localization phenomena will also occur beyond this setting. 
In particular, we expect that similar results should be possible for higher-dimensional wave equations and for control mechanisms acting on lower-dimensional subsets, such as boundary parts or interior interfaces. 
It is to be expected that localized perturbations lead to a localized response if and only if the control geometry prevents waves from travelling arbitrarily far without interacting with the control domain.

\noindent The contribution of this paper is threefold. 
First, we formulate a two-dimensional numerical test problem which transfers the localization question from distributed control of one-dimensional wave equations to boundary and line controls on expanding square domains. 
Second, we propose a finite element and implicit midpoint discretization together with a condensed solution strategy for the resulting linear-quadratic optimal control problem, avoiding the assembly of the full space-time KKT system. 
Third, we perform an extensive simulation study where we compare different control geometries numerically and investigate whether regularly distributed line controls lead to a spatial localization behaviour that is uniform with respect to an increasing domain size.

\noindent The paper is structured as follows: In~\Cref{Sec: TheoreticalBackground} we briefly recall the results from the literature which motivated this work.
In~\Cref{Sec: ModelProblem} we introduce the model problem which is considered in our numerical experiments.
In~\Cref{Sec: Discretization} we discuss the discretization technique which we used, to make the simulations computationally feasible.
Finally we discuss the results of our simulations in~\Cref{Sec: SimulationResults}.
 
\section{Theoretical background and working hypothesis}
\label{Sec: TheoreticalBackground}

As mentioned in~\Cref{Sec: Introduction} our numerical study is motivated by recent results on the spatial localization of perturbations in linear-quadratic optimal control problems~\cite{Oppeneiger2026}. 
In this section we will briefly recall the core ideas of these results and formulate the hypothesis which will be investigated numerically in the remainder of the paper.
We first introduce the abstract optimal control problem considered in that work; see \Cref{Sec: AbstractProblem}. We then summarize the central spatial decay result; see \Cref{Sec: SpatialDecayPrevious}. Finally, we formulate the working hypothesis for the present numerical study in \Cref{Sec: Hypothesis}.

\subsection{Abstract optimal control problem}
\label{Sec: AbstractProblem}
We consider a family of optimal control problems indexed by spatial domains
$\Omega$ in a set of domains $\mathcal{O}$.
Each of these optimal control problems may be written in the form
\begin{flalign}
\tag{$\mathrm{OCP}_\Omega ^\mathrm{T}$}
\label{Eq: AbstractOCP}
    \begin{split}
    \underset{(x,u)}{\min }\,\,\, \frac{1}{2} &\int _0^T \n{C \left(x(t)-x^\mathrm{ref}\right)}_{Z}^2 \, + \,  \alpha \n{u(t)-u^\mathrm{ref}}_{U}^2 \, \mathrm{d}t\\
    \textrm{s.t.}: \quad &\dot{x}(t) = A x(t) + B u(t) + F, \,\, x(0)=x^0,
    \end{split}
\end{flalign}
where
\begin{itemize}
    \item $\Omega \in \mathcal{O}\subset \{\Omega \subset \R^d : \Omega \textrm{ open, bounded with Lipschitz boundary}\}$ is a domain of dimension $d\in \N$
    \item the state space $X$, the input space $U$ and the output space $Z$ are Hilbert spaces,
    \item $A: X \supset D(A)\rightarrow X$ generates a strongly-continuous semigroup $(\mathcal{T}(t))_{t\geq 0}$ on~$X$,
    \item the input and output operators $B \in L(U,X)$ and $C \in L(X,Z)$ are bounded,
    \item $\alpha >0$ is the control weight,
    \item the source term $F \in L^1(0,T;X)$ is integrable and
    \item we have an initial value $y^0 \in X$ and references $x^\mathrm{ref} \in X$ and $u^\mathrm{ref} \in U$.
\end{itemize}
All of the above quantities may depend on the spatial domain $\Omega$, i.e. formally we should write $X_\Omega$ instead of $X$, $A_\Omega$ instead of $A$ etc. However in order to keep the notation simple, we omit the $\Omega$-index in the following except where it is indispensable.
The first-order optimality system corresponding to the optimal control problem~\eqref{Eq: AbstractOCP} is
\begin{equation}
\label{Eq: KKT}
    \begin{pmatrix}
        C^*C & 0 & -\frac{\mathrm{d}}{\mathrm{d}t} - A ^*\\
        0 & 0 & E^T\\
        0 & \alpha & -B^*\\
        \frac{\mathrm{d}}{\mathrm{d}t} - A & -B & 0\\
        E^0 & 0 & 0
    \end{pmatrix}
    \begin{pmatrix}
        x\\ u\\ \lambda 
    \end{pmatrix}
    =
    \begin{pmatrix}
        C^*Cy^\mathrm{ref}\\
        0\\
        \alpha u^\mathrm{ref}\\
        F\\
        x^0
    \end{pmatrix}
    \implies
    \underset{:=\mathcal{M}}{\underbrace{
    \begin{pmatrix}
        C^*C & -\frac{\mathrm{d}}{\mathrm{d}t} - A ^*\\
        0 & E^T\\
        \frac{\mathrm{d}}{\mathrm{d}t} - A & -\alpha ^{-1}BB^*\\
        E^0 & 0
    \end{pmatrix}}}
    \underset{:=\xi}{\underbrace{
    \begin{pmatrix}
        x\\ \lambda 
    \end{pmatrix}
    }}
    =
    \underset{:=h}{\underbrace{\begin{pmatrix}
        C^*Cy^\mathrm{ref}\\
        0\\
        Bu^\mathrm{ref} + F\\
        x^0
    \end{pmatrix}}},
\end{equation}
i.e. if $(x,u)\in C(0,T;X)\times L^2(0,T;U)$ is optimal for~\eqref{Eq: AbstractOCP}, then there exists a Lagrange multiplier 
$\lambda \in C(0,T;X)$ such that the optimality system~\eqref{Eq: KKT} is satisfied.
The operators $E^0$ and $E^T$ encode the initial and terminal conditions, that is,
\begin{equation*}
    E^0: C(0,T;X) \rightarrow X, \platz E^0 x := x(0)
    \qquad
    \mathrm{and}
    \qquad
    E^T: C(0,T;X) \rightarrow X, \platz E^T \lambda := \lambda(T).
\end{equation*}

\begin{example}[Wave equation]
    \label{Ex: WaveEq}
    Consider a wave equation with distributed control and Neumann boundary condition
    \begin{align*}
        \ddot{y} - c^2 \Delta y - \chi_{\Omega _C} u - f&= 0,
        && \mathrm{on} \, (0,T)\times\Omega,\\
        \partial _n y&=0, && \mathrm{on} \, \partial \Omega\\
        y(0, \cdot)=y _0, \, \dot{y}(0,\cdot)&=v _0, && \mathrm{on} \, \Omega
    \end{align*}
    where \(\Omega \subset \mathbb{R}^d\) is a bounded Lipschitz domain and $(x_0,v_0) \in H^1(\Omega)\times L^2(\Omega)$. 
    Introducing the first-order state $x = (y,v) := (y,\dot{y})$ the wave equation may be written in the abstract form $\dot{x} = Ax + Bu + F$ from~\Cref{Eq: AbstractOCP} where we have the source term $F = (0,f)$, the differential operator
    \begin{equation*}
        A: D(A):=H^2(\Omega)\times H^1(\Omega) \subset X := H^1(\Omega) \times L^2(\Omega) \rightarrow X, 
        Ax := \begin{pmatrix}
            0 & I\\
            c^2\Delta & 0
        \end{pmatrix} \begin{pmatrix}
            y\\v
        \end{pmatrix}
    \end{equation*}
    and the input/output operators
    \begin{equation*}
        B: U := L^2(\Omega _C) \rightarrow X, \,\, Bu =
        \begin{pmatrix}
            0\\
            \chi_{\Omega _C}u
        \end{pmatrix} \quad \textrm{and} \quad
        C: X \rightarrow L^2(\Omega _O),\,\,  Cx = v |_{\Omega _O}
    \end{equation*}
    with control domain $\Omega _C$ and observation domain $\Omega _O$.
\end{example}

\subsection{Spatial decay of perturbations}
\label{Sec: SpatialDecayPrevious}
In the following we briefly explain the localization mechanism for perturbations in distributed optimal control of the wave equation described in~\cite{Oppeneiger2026}. This mechanism consists of three steps:
First assuming a regular distribution of the control and observation areas over the whole spatial domain, a \textit{domain-uniform} stabilizability and detectability property may be proven for the wave equation with distributed control/observation.
Second, this domain-uniform stabilizability/detectability implies a uniform (in the size of the spatial domain and the control horizon) bound on the solution operator $\mathcal{M}$ of the associated KKT system~\eqref{Eq: KKT}.
Finally using these domain-uniform bounds it is possible to show the following localization result: Perturbations which are exponentially localized around some point $P$ in space also only have an exponentially localized influence on the solution of the optimal control problem~\eqref{Eq: OCP}. In particular this applies to the optimal state, adjoint state, and control.

\noindent For Step 1 we first need to recall the meaning of domain-uniform stabilizability/detectability as defined in~\cite[Definition 12]{Oppeneiger2026}.

\begin{definition}[Domain-uniform exponential stabilizability/detectability]
\label{Def: DomUniStabDet}
We call the family of differential and input operators $(A_\Omega, B_\Omega)_{\Omega \in \mathcal{O}}$ domain-uniformly stabilizable in $\mathcal{O}$ if and only if there exists a family of uniformly (in $\Omega \in \mathcal{O}$) bounded feedback operators $(K_\Omega^B)_{\Omega \in \mathcal{O}}$, $K_\Omega^B \in L(X_ \Omega,U_\Omega)$, $\Omega\in \mathcal{O}$ such that the following two conditions are satisfied:
\begin{itemize}
    \item[(i)] The closed-loop operators $A_\Omega + B_\Omega K_\Omega^B$, $\Omega\in \mathcal{O}$, generate strongly continuous semigroups $(\mathcal{T}_\Omega(t))_{t\geq 0}$, on $X _\Omega$, $\Omega \in \mathcal{O}$.
    \item[(ii)] There exist constants $M,k>0$ such that 
    $\forall \, \Omega\in \mathcal{O} \, \forall \, t \geq 0: \n{\mathcal{T}_\Omega(t)}_{L(X_\Omega)} \leq Me^{-kt}$.
\end{itemize}
We call $(A_\Omega, C_\Omega)_{\Omega \in \mathcal{O}}$, $C_\Omega \in L(X_\Omega, Z_\Omega )$, $\Omega\in \mathcal{O}$ domain-uniformly exponentially detectable in $\mathcal{O}$ if and only if $(A_\Omega^*, C_\Omega^*)$ is domain-uniformly exponentially stabilizable.
\end{definition}

\noindent In~\cite[Theorem 32]{Oppeneiger2026} domain-uniform stabilizability of the wave equation with distributed control is characterized for one-dimensional spatial domain in $\mathcal{O}_\mathrm{1D} :=\{(0,L): L>0\}$, where the control domain is defined by $\Omega _L^C:= \Omega _L\cap \Omega _C$ for $\Omega _L \in \mathcal{O}$ and a total control domain $\Omega _C \subset \R _{\geq 0}$. The following~\Cref{Thm: StabilizabilityWave} directly follows from this result.

\begin{thm}[Domain-uniform stabilizability of the wave equation]
\label{Thm: StabilizabilityWave}
Assume that there exist constants $c_0, c_1>0$ such that for all intervals $I\subset \R _{\geq 0}$ the inequality
$|\Omega _C \cap I| \geq c_1 |I| - c_2$
is fulfilled and $\forall \, L>0: |\Omega _L^C| = |\Omega _L \cap \Omega _C| >0$.
Then the wave equation with distributed control from~\Cref{Ex: WaveEq} is domain-uniformly stabilizable with regard to $\mathcal{O}_\mathrm{1D}$.
\end{thm}

\noindent If a distributed observation domain is defined by $\Omega _L^O:= \Omega _L\cap \Omega _O$, $\Omega _O \subset \R _{\geq 0}$ then domain-uniform detectability is implied by the analogous condition. This follows directly from duality.
The main idea of~\Cref{Thm: StabilizabilityWave} is the following:
The control region is required to be uniformly distributed in the sense that the distances between neighbouring controlled parts remain uniformly bounded. At the same time, the controlled parts must not degenerate locally: in each prescribed neighbourhood, their size has to be bounded from below. This non-degeneracy is needed because the input operator is bounded, so arbitrarily small control regions cannot provide a uniform control effect.
In the following let $W_\Omega ^2 := ((L_{\mu _\varepsilon} ^2(0,T;X_\Omega)^2 \times (X_\Omega)^2)$.
The second step is performed by using the following~\Cref{Thm: Boundedness} (see~\cite[Theorem 12]{Oppeneiger2026}).
\begin{thm}
\label{Thm: Boundedness}
    Let $(A_\Omega, B_\Omega)_{\Omega \in \mathcal{O}}$ and $(A_\Omega, C_\Omega)_{\Omega \in \mathcal{O}}$ be domain-uniformly stabilizable respectively domain-uniformly detectable families
    Then there exists a constant $c>0$ such that 
    $\forall\, T>0 \, \forall \, \Omega \in \mathcal{O}: \n{\left(\mathcal{M}_\Omega ^T\right)^{-1}}_{L(W_\Omega ^2,W_\Omega ^2)}\leq c$.
\end{thm}

\noindent Now we have all the necessary ingredients to describe the spatial decay of perturbations in optimal control in the third step. 
For this purpose consider the perturbed version $\mathcal{M}x^d = h+\varepsilon$ and the error system $\mathcal{M}\delta \xi^d = \varepsilon$ associated with the optimality system~\eqref{Eq: KKT} where $\varepsilon \in (L^1 (0,T;X)\times L^2(X))^2$ is the perturbation, $\xi^d = (x^d, \lambda ^d)$ is the perturbed solution and $\delta \xi^d := \xi^d - \xi = (\delta x, \delta \lambda)$ is the deviation caused by the perturbation.
We call the family $\left(\varepsilon _\Omega ^T\right) _{\Omega \in \mathcal{O}, T>0}$ exponentially localized around some point $P \in \R^d$ in $W^2$ if and only if there exist constants $\mu _\varepsilon, C_\varepsilon >0$ such that
\begin{equation*}
    \forall \, \Omega \in \mathcal{O}\, \forall \, T>0: \n{\varepsilon _\Omega^T}_{W_{\Omega,\mu} ^2} := \n{e^{\mu _\varepsilon \n{P-\cdot}}\varepsilon _\Omega^T}_{W_\Omega ^2} < C_\varepsilon < \infty.
\end{equation*}

\noindent The following~\Cref{Thm: SpatialDecay} is a consequence of~\cite[Theorem 10]{Oppeneiger2026}. It holds true under the assumption, that the differential operator $A$ behaves similar to a first-order differential operator (see~\cite[Assumption 6]{Oppeneiger2026}). For the purpose of this work it suffices to say, that the wave equation satisfies this Assumption (see~\cite[Example 8]{Oppeneiger2026}).
\begin{thm}[Spatial decay of local perturbations]
\label{Thm: SpatialDecay}
    Assume that there exists a constant $c>0$ such that
    \begin{equation}
    \label{eq: SolOperatorBoundedness}
        \forall \, T>0 \, \forall \, \Omega \in \mathcal{O}: \n{\left(\mathcal{M} _\Omega ^T \right) ^{-1}}_{L(W_\Omega ^2, W_ \Omega ^2)}\leq c.
    \end{equation}
    Let the family of perturbations $\left(\varepsilon _\Omega ^T\right) _{\Omega \in \mathcal{O}, T>0}$ be exponentially localized around $P \in \R^d$ in $W^2$ with constants $\mu _\varepsilon, C_\varepsilon >0$
    Then, there exist constants $K, \mu >0$ such that for all $T>0$ and for all $\Omega \in \mathcal{O}$ the inequality
    \begin{equation}
        \label{SpatialDecayTheorem}
        \n{e^{\mu \n{P-\cdot}_1} \delta x_\Omega^T}_{L^2(0,T;X)} + 
        \n{e^{\mu \n{P-\cdot}_1} \delta \lambda_\Omega^T}_{L^2(0,T;X)} + 
        \n{e^{\mu \n{P-\cdot}_1} \delta u_\Omega^T}_{L^2(0,T;U)} 
        \leq K \n{e^{\mu \n{P-\cdot}_1} \varepsilon_{\textcolor{black}{\Omega}}^T}_{W^2} 
        \leq K\cdot C_\varepsilon,
    \end{equation}
    i.e. the deviation in the optimal state, optimal control and adjoint state is exponentially localized around $P$ in $L^2$.
\end{thm}

\subsection{working hypothesis for the numerical study}
\label{Sec: Hypothesis}
The results presented in~\Cref{Sec: SpatialDecayPrevious} come with two main drawbacks. 
First the sufficient condition for domain-uniform stabilizability and detectability in~\Cref{Thm: StabilizabilityWave} only applies to one-dimensional domains. 
Second the results on the domain-uniform boundedness of the solution operator in~\Cref{Thm: Boundedness} as well as on the spatial decay of perturbations in~\Cref{Thm: SpatialDecay} were only rigorously shown for bounded input and output operators. 
For hyperbolic equations like the wave equation this disqualifies boundary control and observation.

\noindent The working hypothesis for the present paper is that similar results may also be shown for scenarios with higher-dimensional domains and/or boundary control observation.
For a wave equation one could suspect that - similar to the geometric control condition (GCC)~\cite{Bardos1992} - it may suffice in such a scenario to ensure, that the distance a wave can travel in a straight line without touching any control hyperplanes is bounded from above.

\noindent In order to test this hypothesis we have implemented a wave equation on a square domain with a rectangular grid of line controls.

\noindent If our hypothesis is correct, we should observe a spatial decay of perturbations under the influence of an optimal control on the grid lines, which is independent of the domain size.
By contrast, if the control only acts on the outer boundary of the square domain, then the spatial decay is expected to deteriorate as the domain size grows.
This seems plausible since a perturbation must propagate across the entire spatial domain before it can be damped by the optimal control.

\section{Model problem}
\label{Sec: ModelProblem}
\begin{figure}[h!]
    \centering
    \includegraphics[width=0.5\linewidth]{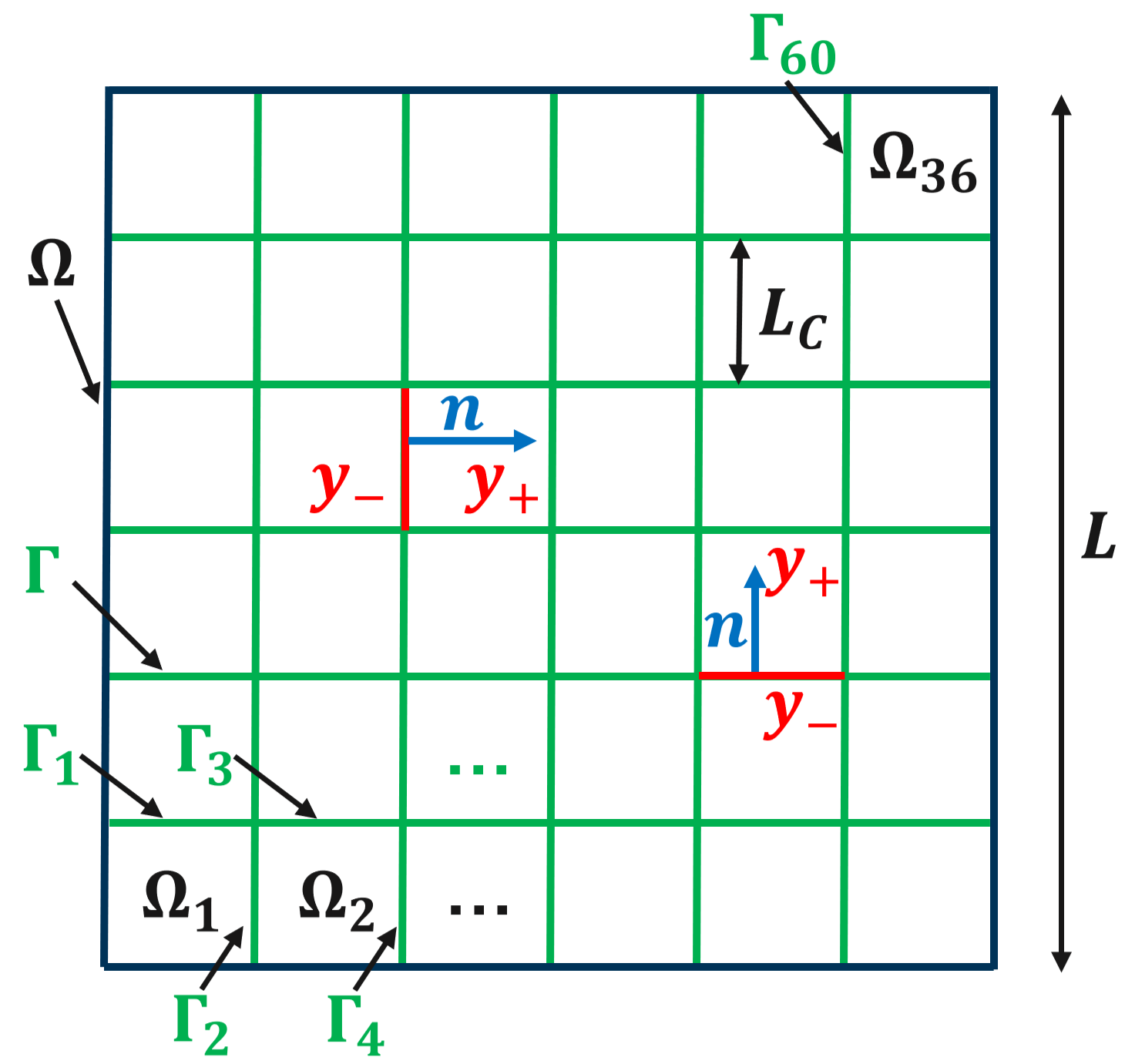}
    \caption{Visualization of the spatial domain $\Omega$, the subdomains $\Omega _k$ and the grid $\Gamma$ for $N=6$.}
    \label{fig: SpatialDomain}
\end{figure}

For our simulation study we consider an optimal control problem for the wave equation given by
\begin{subequations}
\label{Eq: OCP}
\begin{align}
\min_{y,u} \,
  \frac{1}{2}\int _0^T \!\! \|y(\cdot, t)-y_d\|_{H}^2 &+ k_v \|\dot{y}(\cdot, t)-v_d\|_{H}^2 + \|u(\cdot,t)\|_U^2
\label{Eq: Objective}
\\
\text{s.t.} \quad
 \ddot{y} - c^2 \Delta y &= 0   && \text{on } \left(\Omega \setminus \Gamma\right)\times [0,T],
\label{Eq: State}
\\
\label{Eq: Neumann}
\partial _n y &= 0   && \text{on } \partial \Omega \setminus \Gamma _C,\\
\label{Eq: Init}
y(\cdot, 0) = y_0, v(\cdot, 0) &= v_0   && \text{on } \Omega \setminus \Gamma,\\
\label{Eq: Continuity}
y_+- y_- &= 0   && \text{on } \Gamma,\\
\label{Eq: BoundaryControl}
\partial _n y_+ - \partial _n y_- &= u   && \text{on } \Gamma _C,
\end{align}
\end{subequations}
where
\begin{itemize}
    \item $\Omega = (0,L)^2 \subset \R^2, L = NL_C$ is a square domain partitioned into subsquares as depicted in~\Cref{fig: SpatialDomain},
    \item $\Gamma \subset \Omega$ denotes the inner boundaries of the subsquares $\Omega _1, \ldots, \Omega _{N^2}$, i.e.$\Gamma = \left(\cup_{k=1}^{N^2} \partial \Omega _k \right) \setminus \partial \Omega$,
    \item $\Gamma _C \subset \partial \Omega \cup \Gamma$ is the control domain on which the boundary/line control is acting,
    \item $H = L^2(\Omega)$ and $U = L^2(\Gamma _C)$ are the spaces of square integrable functions and inputs,
    \item $n = n(\omega) \in \R^2$ denotes the normal vector in a point $\omega \in \Gamma$ defined as in~\Cref{fig: SpatialDomain} and $c \in \R$ is the propagation velocity,
    \item $y_0 \in H^1(\Omega)$ and $v_0 \in L^2(\Omega)$ denote the initial displacement and velocity distributions on the spatial domain $\Omega$,
    \item $y_d \in X$ and $v_d \in H$ denote the desired displacement and velocity distributions on the spatial domain $\Omega$,
    \item $y_+$ and $y_-$ are the solution along a boundary in and against the direction of the normal vector $n$, i.e. we have 
    \begin{equation*}
        \forall \, \omega \in \Gamma: y_+(\omega) = \lim_{\varepsilon \downarrow 0} y(\omega + \varepsilon n) \quad \mathrm{and} \quad y_-(\omega) = \lim_{\varepsilon \downarrow 0} y(\omega - \varepsilon n),
    \end{equation*}
    \item and $\partial _n y_+ (\omega) := \left(\nabla y_+\right)^\top n$ and $\partial _n y_- (\omega) := \left(\nabla y_-\right)^\top n$ are the normal derivatives of the displacement $y$ on both sides of the inner and outer boundary facets.
\end{itemize}
For $\omega \in \partial \Omega$ we define $\partial _n y_- (\omega) = 0$ if $n$ points into $\Omega$ and $\partial _n y_+ (\omega) = 0$ if $n$ points out of $\Omega$.

\noindent To investigate the localization properties of this optimal control we will consider a perturbation of the initial value, which is given by a gaussian impulse of the form
\begin{equation}
\label{Eq: Perturbation}
    d(\omega) = e^{- \frac{1}{2 \sigma ^2} \n{P-\omega}_2^2}
\end{equation}
where $P \in \R^2$ is the mean value around which the perturbation is localized and $\sigma >0$ may be seen as a standard deviation.

\noindent Furthermore we will consider two different control scenarios: In the first the control only acts on the outer boundary of the square domain, i.e. $\Gamma _C = \partial \Omega$. In the second scenario the control acts on the outer boundary as well as on the inner control grid defined by the boundaries of the subdomains, i.e. $\Gamma _C = \partial \Omega \cup \Gamma$.

\noindent To confirm our working hypothesis, domain-uniform localization should only be visible in the second scenario, since the distance between the control lines grows with the domain size in the first.

\section{Discretization of the model problem}
\label{Sec: Discretization}
We now describe the discretization used for the numerical solution of the model problem introduced above. The spatial discretization is based on conforming first-order finite elements on a structured triangular mesh, while time integration is performed by the implicit midpoint rule. This leads to a fully discrete linear-quadratic optimal control problem, which is subsequently solved in condensed form rather than by assembling the full space-time KKT system.

\subsection{Discretization of the state equation}
\label{Sec: StateDiscretization}
In order to discretize the model problem we first derive its weak formulation using Greens formula and the definition of the normal vectors depicted in~\Cref{fig: SpatialDomain}. For $y \in H^2(\Omega), w \in H^1(\Omega)$ we find
\begin{align*}
    c^2 \left\langle \Delta y, w \right\rangle_{L^2(\Omega)}
    &= \sum_{k=1}^{N^2} c^2 \left\langle \Delta y, w \right\rangle_{L^2(\Omega_k)}
    = \sum_{k=1}^{N^2} c^2 \left( \left\langle \partial_{n_\mathrm{out}} y, w \right\rangle_{L^2(\partial\Omega_k)} - \left\langle \nabla y, \nabla w \right\rangle_{L^2(\Omega_k)} \right)\\
    &= c^2 \left( \left\langle \partial_n y, w \right\rangle_{L^2(\partial\Omega)} +
        \sum_{k=1}^{2N(N-1)} \left\langle \partial_n y_{-} - \partial_n y_{+}, w \right\rangle_{L^2(\Gamma _k)} - \sum_{k=1}^{N^2} \left\langle \nabla y, \nabla w \right\rangle_{L^2(\Omega_k)} \right) \\
    &= - c^2 \left( \left\langle u, w \right\rangle_{L^2(\Gamma _C)} + \left\langle \nabla y, \nabla w \right\rangle_{L^2(\Omega)} \right),
\end{align*}
where $n_\mathrm{out}$ is an outer normal vector with regard of $\Omega _k$ and the last step follows from equations~\eqref{Eq: Neumann} and~\eqref{Eq: BoundaryControl}. Therefore the weak version of the above wave equation with boundary/line control is
\begin{equation}
    \label{WeakFormulation}
    \forall \,  w \in H^1 (\Omega): \left\langle \ddot{y}, w \right\rangle_{L^2(\Omega)} + c^2 \left\langle u, w \right\rangle_{L^2(\Gamma _C)} + c^2 \left\langle \nabla y, \nabla w \right\rangle_{L^2(\Omega)}
    = 0.
\end{equation}
We discretize the square domain $\Omega_L=(0,L)^2$ by a structured triangular mesh. Starting from a uniform Cartesian grid with mesh size $h$, each rectangular cell is subdivided into two triangles. On this triangulation, we use conforming finite elements of first order, i.e. continuous, piecewise affine $P_1$ basis functions $\varphi _1, \ldots , \varphi _{n_m}$ leading to a mesh with $n_\mathrm{m}$ degrees of freedom. We denote the set of nodes by $p_h$ and the set of triangular elements by $t_h$.
Thus we arrive at a second-order semidiscrete version the wave equation, which is
\begin{equation*}
    w_h^\top M_h \ddot{y}_h + c^2 w_h^\top K_h x_h + c^2 w_h^\top \tilde{B} _h u_h = 0,
\end{equation*}
where $M_h, K_h \in \R^{n_m \times n_m}$ denote the standard mass and stiffness matrices corresponding to the bilinear forms $m(w,v):=\scp{w}{v}_{L^2(\Omega)}$ and $a(w,v):=\scp{\nabla w}{\nabla v}_{L^2(\Omega)}$. The discretized control matrix $\tilde{B}_h$ is defined by $\tilde{B}_h := \Sigma _h(:,I_C)$ where $I_C$ denotes the indices of the nodes corresponding to the control domain and $\Sigma _h \in \R^{n_m \times n_m}$ is the mass matrix with regard to the control domain, i.e. $\Sigma _h (i,j):= \scp{\varphi _i}{\varphi _j}_{L^2(\Gamma _C)}$.
Introducing a new state $x_h = (y_h, v_h) := (y_h, \dot{y}_h)$ in analogy to the first-order formulation of the wave equation in~\Cref{Ex: WaveEq} leads to
\begin{equation}
\label{Eq: SemidiscreteWave}
    \underset{=:E_h}{\underbrace{\begin{pmatrix}
        I & 0\\
        0 & M_h
    \end{pmatrix}}}
    \underset{=:\dot{x}_h}{\underbrace{\begin{pmatrix}
        \dot{y}_h\\
        \dot{v}_h
    \end{pmatrix}}}
    = \underset{=:A_h}{\underbrace{\begin{pmatrix}
        0 & I\\
        -K_h & 0
    \end{pmatrix}}}
    \underset{=:x_h}{\underbrace{\begin{pmatrix}
        y_h\\
        v_h
    \end{pmatrix}}}
    \underset{=:B_h}{\underbrace{\begin{pmatrix}
        0\\
        \tilde{B}_h
    \end{pmatrix}}}
    u_h, \qquad y_h(0) \approx y_0, \quad v_h(0) \approx v_0
\end{equation}
\noindent For the time discretization we use the implicit midpoint rule. This choice is
motivated by its structure-preserving properties: the implicit midpoint rule is symplectic and preserves quadratic invariants (i.e. quantities constant in time) exactly.
Since the semi-discrete, uncontrolled wave equation is a linear Hamiltonian
system with a quadratic energy, the method preserves the discrete wave energy
in the absence of forcing and control; see, e.g.,
\cite{Hairer2006,Ascher1999}. This way we are able to avoid artificial numerical dissipation.
For the semi-discrete system in~\eqref{Eq: SemidiscreteWave}, sampling time $\Delta t$ and time steps $t_0=0,\ldots ,t_{n_t} = n _t \Delta t$ the implicit midpoint rule reads
\begin{equation}
    E_h \frac{x_h^{n+1} - x_h^n}{\Delta t} = A_h \frac{x_h^{n+1} - x_h^n}{2} + B_hu_h^n \implies \underset{=:D_+}{\underbrace{\left( \frac{1}{\Delta t} E_h - \frac{1}{2}A_h \right)}}x_h^{n+1} = \underset{=:-D_-}{\underbrace{\left( \frac{1}{\Delta t} E_h + \frac{1}{2}A_h \right)}}x_h^{n} + B_h u_h^n.
\end{equation}
Combining all time steps we find the linear equation system
\begin{equation}
\label{Eq: StateDiscrete}
    \underset{=:\hat{A}_h}{\underbrace{\begin{pmatrix}
        I & 0 & \hdots & 0\\
        D_- & D_+ & \hdots & 0\\
        & \ddots & \ddots & \vdots\\
        && D_- & D_+
    \end{pmatrix}}}
    \underset{=:X_h}{\underbrace{\begin{pmatrix}
        x_h^0\\
        \vdots \\
        \vdots \\
        x_h^{n_t}
    \end{pmatrix}}}
    +
    \underset{=:\hat{B}_h}{\underbrace{\begin{pmatrix}
        0 &  \hdots & 0\\
        -B_h &  &  \\
        &   \ddots & \\
        &&   -B_h
    \end{pmatrix}}}
    \underset{=:U_h}{\underbrace{\begin{pmatrix}
        u_h^0\\
        \vdots \\
        \vdots \\
        u_h^{n_t-1}
    \end{pmatrix}}}
    =
    \underset{=:F_h}{\underbrace{\begin{pmatrix}
        \begin{pmatrix}
            y_0\\
            v_0
        \end{pmatrix}\\
        0 \\
        \vdots \\
        0
    \end{pmatrix}}}
\end{equation}

\subsection{Discretization of the cost functional}
\label{Sec: CostDiscretization}

In analogy to~\Cref{Sec: StateDiscretization} the cost functional may be discretized in space using finite elements. The semidiscrete cost is
\begin{equation*}
    J_h(x_h,u_h):=\frac{1}{2} \int _0^T \n{C_h x_h(t)}_2^2 + k_u \n{u_h}_{R_h}^2 \mathrm{d}t
\end{equation*}
where
\begin{equation*}
    C_h := \begin{pmatrix}
        M_h & 0\\
        0 & k_v M_h
    \end{pmatrix}
    \quad \mathrm{and} \quad
    R_h := \Sigma _h (I_C,I_C).
\end{equation*}
Using the midpoint integration rule for discretixation in time we find the discrete cost functional
\begin{flalign*}
    J_{h,\Delta t}(X_h,U) &:= \frac{\Delta t}{2} \sum _{n=0}^{n_t-1} \left(\frac{x_h^{n+1} + x_h^n}{2}\right)^\top \underset{=:Q_h}{\underbrace{C_h^\top C_h}} \frac{x_h^{n+1} + x_h^n}{2} + \left(u_h^k\right)^\top R_h u_h^k = \frac{1}{2}\left(X_h^\top \hat{Q}_hX_h + U_h^\top \hat{R}_hU_h\right)
\end{flalign*}
where
\begin{flalign*}
    \hat{Q}_h:= \frac{\Delta t}{4} 
    \begin{pmatrix}
        Q_h & Q_h &&&\\
        Q_h & 2Q_h &&&\\
        &\ddots&\ddots &\ddots&\\
        &&&2Q_h & Q_h \\
        &&&Q_h & Q_h
    \end{pmatrix}
    \quad \mathrm{and} \quad 
    \hat{R}_h := \frac{\Delta t}{2} k_u
    \begin{pmatrix}
        R_h &  &\\
        &\ddots &\\
        & & R_h
    \end{pmatrix}
\end{flalign*}

\subsection{Discrete optimality system}
Instead of assembling and solving a discretized version of the full space-time KKT system in~\eqref{Eq: KKT}, we use a
condensed formulation in which the state and adjoint equations are eliminated.
This leads to a reduced optimality system posed only in terms of the control.
The action of the reduced Hessian is evaluated matrix-free by one forward and
one adjoint solve.

\noindent This approach has several advantages. First, the large indefinite KKT matrix 
does not have to be assembled or factorized. Second, the reduced
system is symmetric positive definite under the standard assumptions of the
linear-quadratic problem, so that it can be solved by a preconditioned conjugate
gradient (PCG) method rather than by a general Krylov method for indefinite systems like for example MINRES or GMRES.
Third, the time-stepping structure of the wave equation can be exploited directly in each iteration of PCG, since a forward-backward solve only requires the repeated application of $D_+^{-1}$ and $D_+^{-\top}$ to the current state/adjoint state.
This procedure is feasible, since the matrix $\hat{A}$ from the discrete state equation is invertible, which allows us to eliminate the discrete $X_h$ by solving~\eqref{Eq: StateDiscrete} for $U_h$. Thus we find
\begin{flalign*}
    X_h = \hat{A}_h^{-1} F_h - \hat{A}_h^{-1} \hat{B}_h U_h \implies 
    J_{h,\Delta t}(X_h,U_h) &= \frac{1}{2} U_h^\top \left(\hat{B}_h^\top \hat{A}_h^{-\top} \hat{Q}_h \hat{A}_h^{-1} \hat{B}_h + \hat{R}_h\right)U_h + \frac{1}{2}F_h^\top \hat{A}_h^{-\top} \hat{Q}_h \hat{A}_h^{-1} \hat{B}_hU_h.
\end{flalign*}
Hence the discretized optimization problem can be written as
\begin{equation*}
    \underset{U_h \in \R^{n_t \cdot n_C}}{\min} \n{\left(\hat{C}_h\hat{A}_h^{-1}\hat{B}_h + \hat{R}_h^{\frac{1}{2}}\right)U_h}_2^2 + \left(\hat{C}_h\hat{A}_h^{-1}F_h\right)^\top \hat{C}_h\hat{A}_h^{-1} \hat{B}_h U_h
\end{equation*}
Taking the derivative with regard to $U_h$ we arrive at the corresponding necessary optimality condition
\begin{equation}
\label{Eq: CondensedOptimalityCondition}
    \underset{=:\mathfrak{H}_h}{\underbrace{\left(\hat{C}_h\hat{A}_h^{-1}\hat{B}_h + \hat{R}_h^{\frac{1}{2}}\right)^\top \left(\hat{C}_h\hat{A}_h^{-1}\hat{B}_h + \hat{R}_h^{\frac{1}{2}}\right)}} U_h = -\underset{=:\mathfrak{b}_h}{\underbrace{\frac{1}{2}\left(\hat{C}_h\hat{A}_h^{-1}\hat{B}_h\right)^\top\hat{C}_h\hat{A}_h^{-1}F_h}}
\end{equation}
To solve the discrete optimality system~\eqref{Eq: CondensedOptimalityCondition} we first note that applying $\hat{C}\hat{A}^{-1}$ to $F$ is the same as sequentially forward-solving the semidiscrete system~\eqref{Eq: SemidiscreteWave} with initial values $x_0, v_0$ using the implicit midpoint rule. 
In analogy applying $\left(\hat{C}_h\hat{A}_h^{-1}\hat{B}_h\right)^\top$ to the resulting state trajectory $X_0 = \hat{C}_h\hat{A}_h^{-1}F_h$ is the same as solving the adjoint equation of~\eqref{Eq: SemidiscreteWave} backward in time using the implicit midpoint rule. From the perspective of the adjoint equation the state trajectory $X_0$ is an input.
Thus we may obtain the right hand side of~\eqref{Eq: CondensedOptimalityCondition}.
The main computational challenge of this approach lies in the repeated application of $D_+^{-1}$ and $D_+^{-\top}$ to the current state/adjoint state.
To simplify this computation we use the Schur complement. We find
\begin{equation*}
    D_+ y = \begin{pmatrix}
        \frac{1}{\Delta t} I & -\frac{1}{2}I\\
        \frac{1}{2} K_h & \frac{1}{\Delta t} M_h
    \end{pmatrix}
    \begin{pmatrix}
        y_1 \\ y_2
    \end{pmatrix}
    = b = \begin{pmatrix}
        b_1\\b_2
    \end{pmatrix} \Lra y_2 = \frac{2}{\Delta t} y_1 - 2b_1 \land y_1 = S^{-1} \left( \frac{\Delta t}{2} b_2 + M_h b_1\right)
\end{equation*}
where $S:= \frac{1}{\Delta t}M_h + \frac{\Delta}{4} K_h$. In analogy we find
\begin{equation*}
    D_+^\top y 
    = b \Lra y_1 = \Delta t b_1 - \frac{\Delta t}{2} K_h y_2 \land y_2 = S^{-1} \left( \frac{\Delta t}{2} b_1 + b_2\right)
\end{equation*}
Therefore it suffices to factorize the Schur matrix $S$ of size $n_m$ instead of the matrix $D_+$ of size $2n_m$. This factorization only needs to be carried out once. 
Afterwards we can use it to solve the linear equation systems originating from the implicit midpoint rule in each time step by forward backward substitution.
To solve the optimality condition~\Cref{Eq: CondensedOptimalityCondition} we use the preconditioned conjugate gradient method. Here the control weighting matrix $\hat{R}$ serves as preconditioner in order to avoid unfavourable conditioning which may be caused by a very small control weight $k_u$.
The computationally most challenging part of PCG is the application of the reduced hessian matrix $\mathfrak{H}_h$ to a control sequence.
However using the above ideas this simply corresponds to performing a forward-backward solve of the wave equation using the implicit midpoint again.
Therefore the only matrix we need to factorize in order to solve the complete optimal control problem numerically is the Schur matrix $S$.

\section{Simulation results}
\label{Sec: SimulationResults}
In this section, the numerical results obtained for the model problem introduced in~\Cref{Sec: ModelProblem} are presented.
Our main aim is to investigate whether our working hypothesis from~\Cref{Sec: Hypothesis} is reflected in the numerical results.
As described in~\Cref{Sec: ModelProblem}, we compare two different control scenarios:
\begin{itemize}
    \item Scenario 1: The control acts only on the outer boundary of the spatial domain $\Omega$.
    \item Scenario 2: The control acts on the outer boundary as well as on the rectangular grid depicted in~\Cref{fig: SpatialDomain}.
\end{itemize}

\noindent All simulations were carried out with zero reference trajectories, i.e. $y^{\mathrm{ref}}=0$ and $v^{\mathrm{ref}}=0$. 
This way the initial displacement $y_0=d$ defined in~\eqref{Eq: Perturbation} can be interpreted as a local in space perturbation of the zero solution of the wave equation.

\noindent For the numerical experiments the perturbation was centered at $P=(0.65,0.65)$ while the standard deviation was set to $\sigma=0.12$.
Furthermore the velocity was initialized with zero. 
The optimal control problems were solved on a sequence of squares of increasing size $L=0.5,1.0,\ldots,10.0$, with fixed distance $L_C=1$ between neighbouring control lines. 
Since the numerical experiments are intended to provide insight into the continuous optimal control problem, the computations were carried out for several different spatial resolutions $h\in\{0.025, 0.05, 0.1, 0.2, 0.4\}$. 
This allows us to assess whether the observed localization behaviour is a property of the underlying optimally controlled wave equation rather than an artefact of a particular discretization.
\noindent The remaining parameters were kept fixed throughout all simulations with a wave speed of $c=1.0$ and weights $k_u=10^{-3}$ and $k_v=10^{-3}$ in the cost functional.
The final time was set to $T=10.0$ to ensure, that the perturbation is able to travel through the whole spatial domain during each simulation run.
For the time discretization $n_t=150$ time steps were used. 
The reduced optimality system is solved by the preconditioned conjugate gradient method with relative tolerance $10^{-6}$ and a maximum number of $3000$ iterations.

\begin{figure}[h!]
    \centering
    \includegraphics[width=\linewidth]{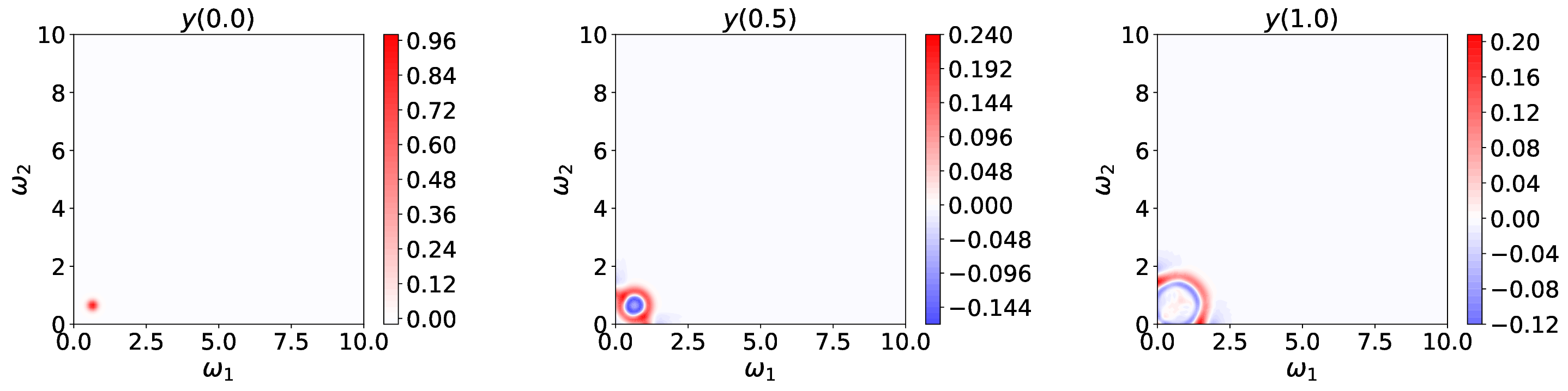}
    \centering
    \includegraphics[width=\linewidth]{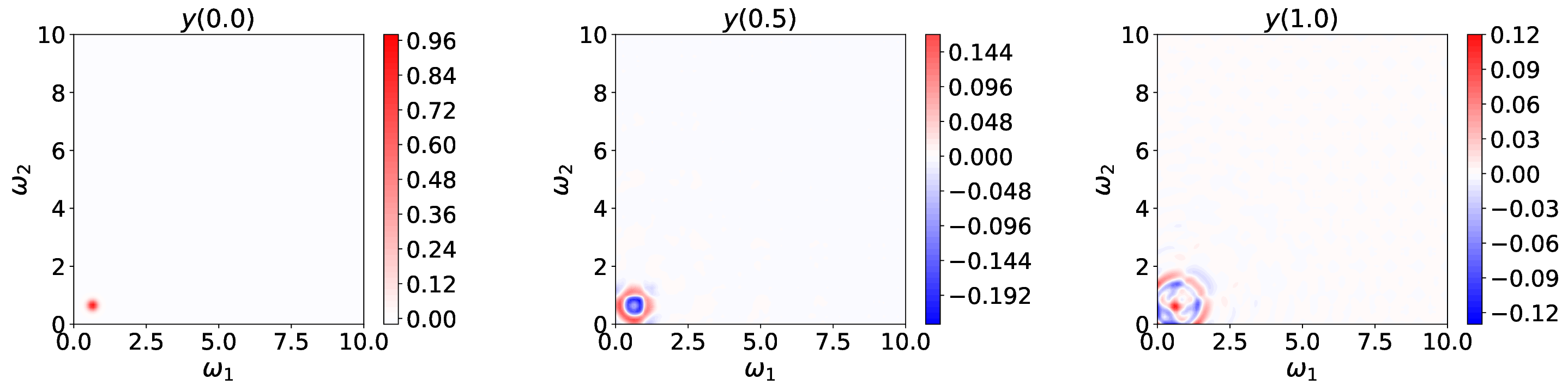}
    \centering
    \includegraphics[width=\linewidth]{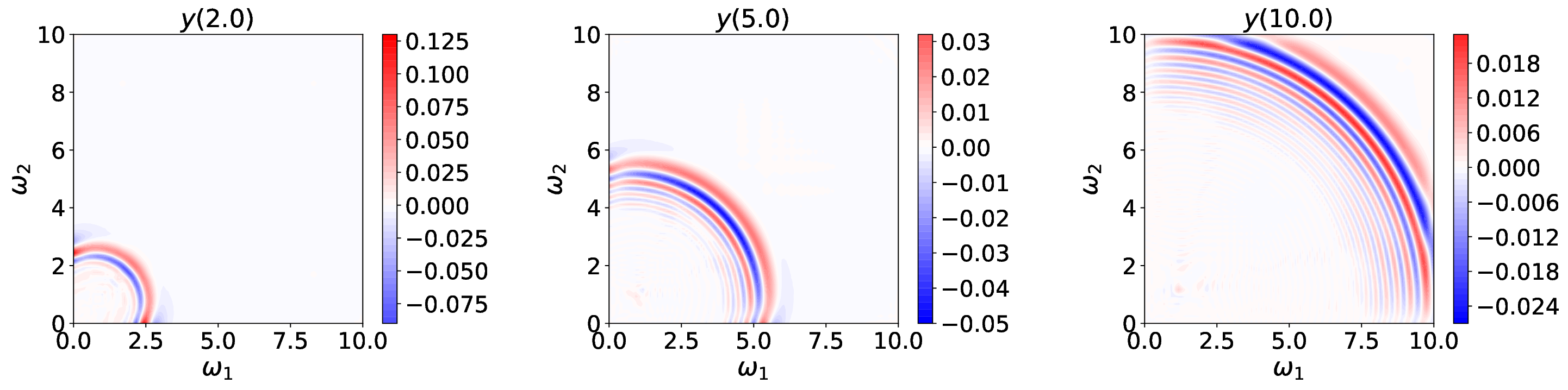}
    \centering
    \includegraphics[width=\linewidth]{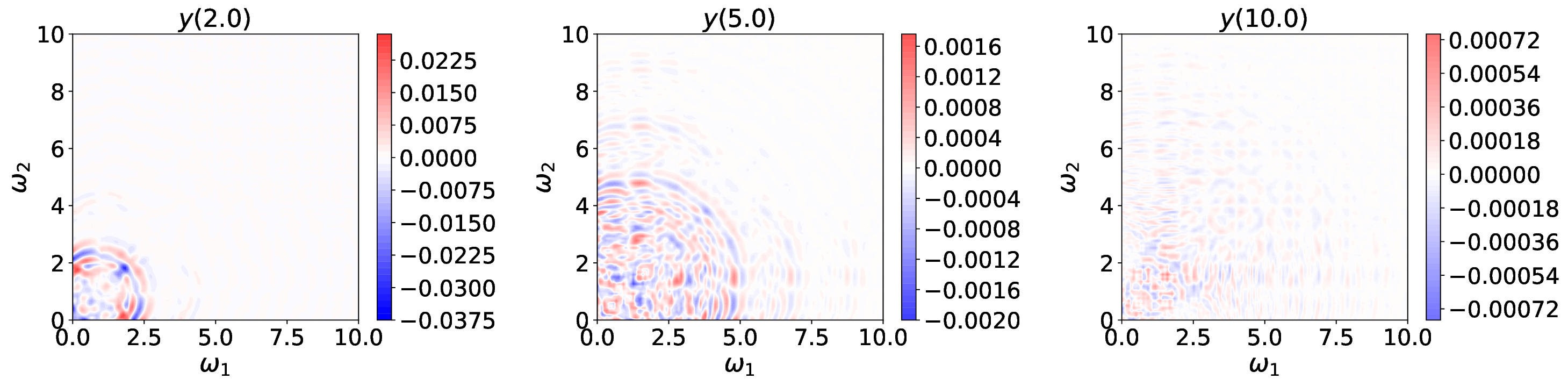}
    \caption{Snapshots of the displacement for pure boundary (top) and grid control (bottom) with resolution $h=0.025$ and domain size $L=10$.}
    \label{Fig: Displacement}
\end{figure}

\noindent \Cref{Fig: Displacement} shows snapshots of the displacement for the optimally controlled wave equation for both scenarios.
This simulation was conducted on large spatial domain ($L=10$) with a fine resolution of $h=0.025$.
If the control acts only at the outer boundary, the initial perturbation is able to propagate over the entire spatial domain as an undamped wave with time.
Since the area over which this wave spreads increases quadratically with time, we observe a slow polynomial decay of the perturbation in space caused by the inherent nature of the problem.
In the scenario, where the control acts on a regularly distributed grid, the decay becomes significantly faster.
Since it is not necessary for the perturbation to spread to the far boundary before it can be damped this is only reasonable.
Instead the wave is cancelled as soon as it reaches the control grid on the boundary of the bottom left subdomain $\Omega _1$.
As a result, at $t=5.0$ the displacement is already so small that it falls within the numerical error tolerance corresponding to the chosen resolution $h=0.025$.

\begin{figure}[htbp]
    \centering
    \includegraphics[width=\linewidth]{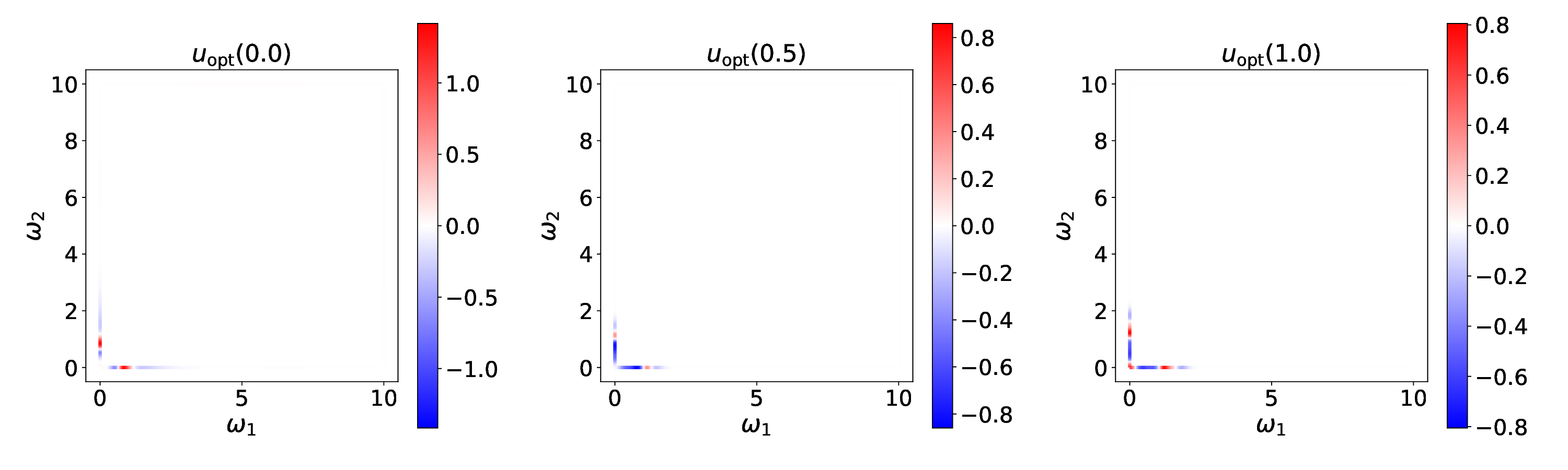}
    \centering
    \includegraphics[width=\linewidth]{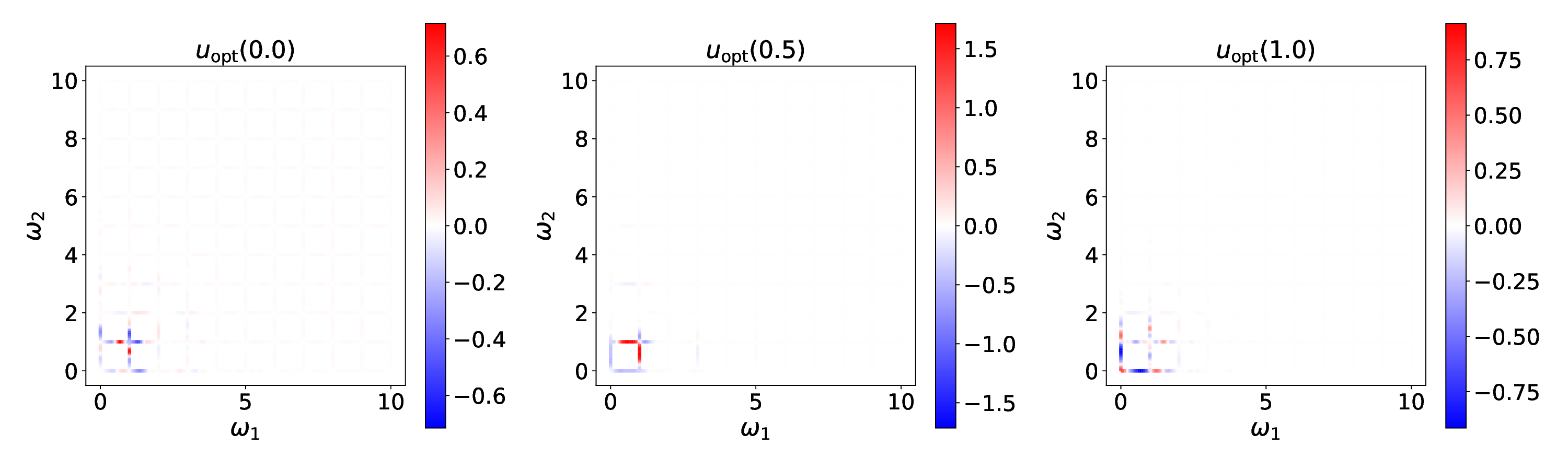}
    \centering
    \includegraphics[width=\linewidth]{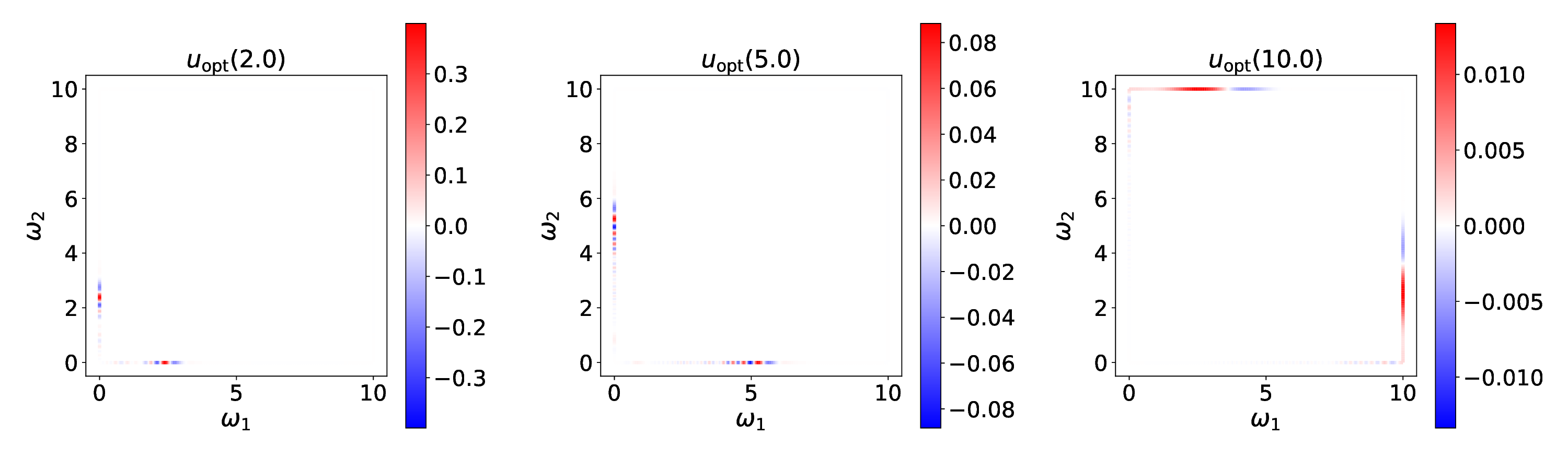}
    \centering
    \includegraphics[width=\linewidth]{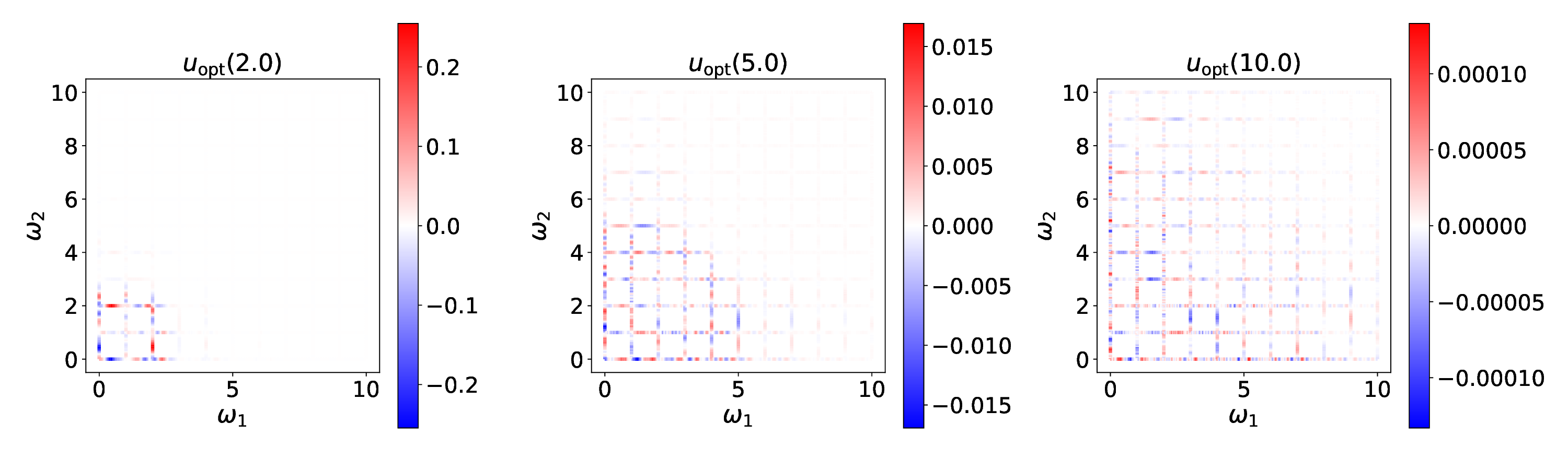}
    \caption{Snapshots of the optimal control for pure boundary (top) and grid control (bottom) with resolution $h=0.025$ and domain size $L=10$.}
    \label{Fig: OptimalControl}
\end{figure}

\noindent The snapshots of the corresponding optimal control in~\Cref{Fig: OptimalControl} show a consistent behaviour.
In Scenario 1 the control is activated in those segments where the wave originating from the initial perturbation reaches the outer boundary of the spatial domain.
In Scenario 2 an analogue phenomenon can be observed for the optimal control on the grid lines which are near the waves current location.

\noindent Once again, we notice that the magnitude of the optimal control declines very rapidly in Scenario 2 due to the rapid damping of the disturbance.
In Scenario 1, this decay occurs much more slowly, and the amplitude of the optimal control travels along the boundary away from the origin.

\begin{figure}[htbp]
    \centering
    \includegraphics[width=0.48\textwidth]{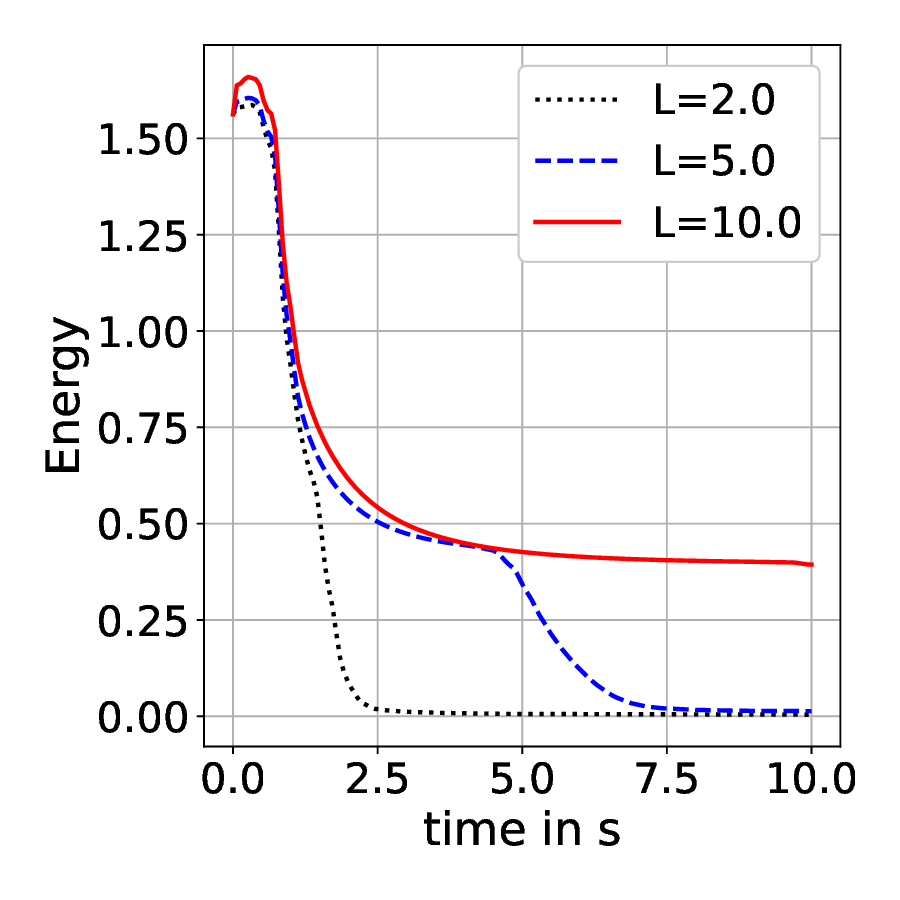}
    \hfill
    \includegraphics[width=0.48\textwidth]{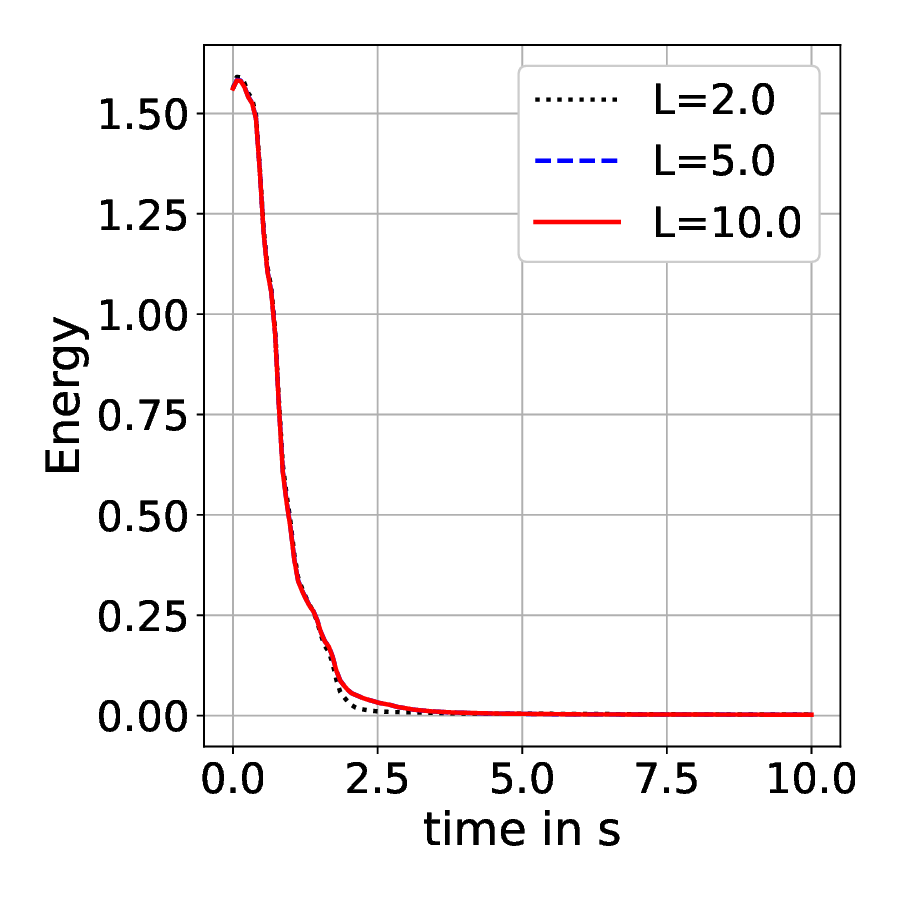}
    \caption{Development of the wave energy over time for different sizes of the spatial domain with pure boundary (left) and grid control (right).}
    \label{Fig: Energy}
\end{figure}
\noindent The difference in damping between the two different scenarios becomes also visible in~\Cref{Fig: Energy} where the wave energy
\begin{equation*}
    E_\mathrm{opt}(t) := \frac{1}{2} \left(\n{\nabla y_\mathrm{opt}(\cdot,t)} _{L^2(\Omega)}
    + \n{v_\mathrm{opt}(\cdot,t)} _{L^2(\Omega)}\right).
\end{equation*}
is depicted for three different sizes of the spatial domain.
If the control acts on the whole grid, we observe a rapid decay of the wave energy which shows hardly any change at all, as the domain size increases.
However if the control only acts on the outer boundary, a significant change of behaviour becomes visible.
For $L=2.0$ the decay in energy is nearly identical to what we found in Scenario 2.
For $L=5.0$ the energy declines in two steps around $t=1$ and $t=5$.
And for $L=10.0$ the prescribed time horizon is not sufficient for the energy to be completely drawn from the system.
This can be explained by the fact that, within this time frame, the wave does not cover the entire diagonal of the spatial domain preventing the control from absorbing the remaining energy.

\begin{figure}[htbp]
    \centering
    \includegraphics[width=0.48\textwidth]{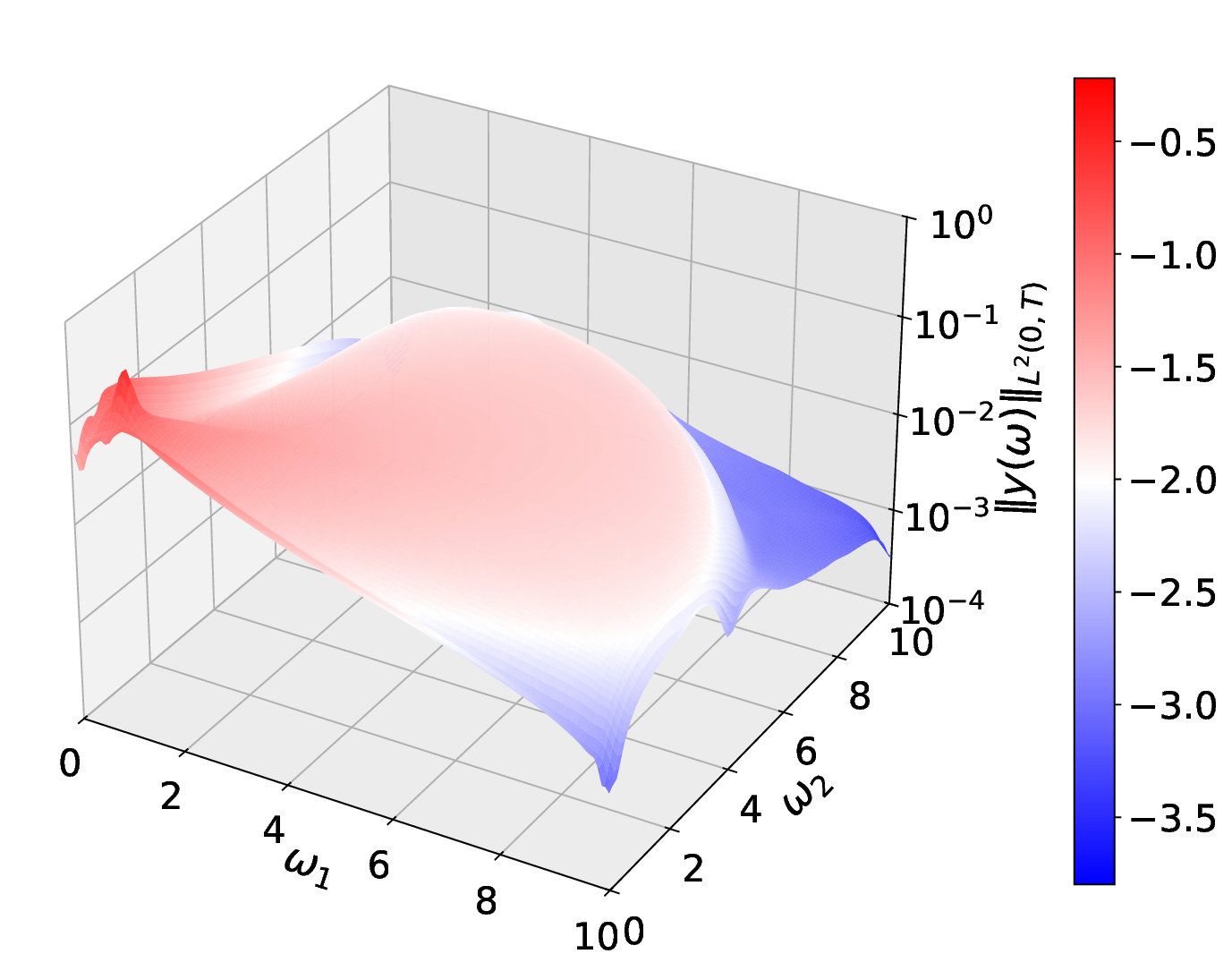}
    \hfill         
\includegraphics[width=0.48\textwidth]{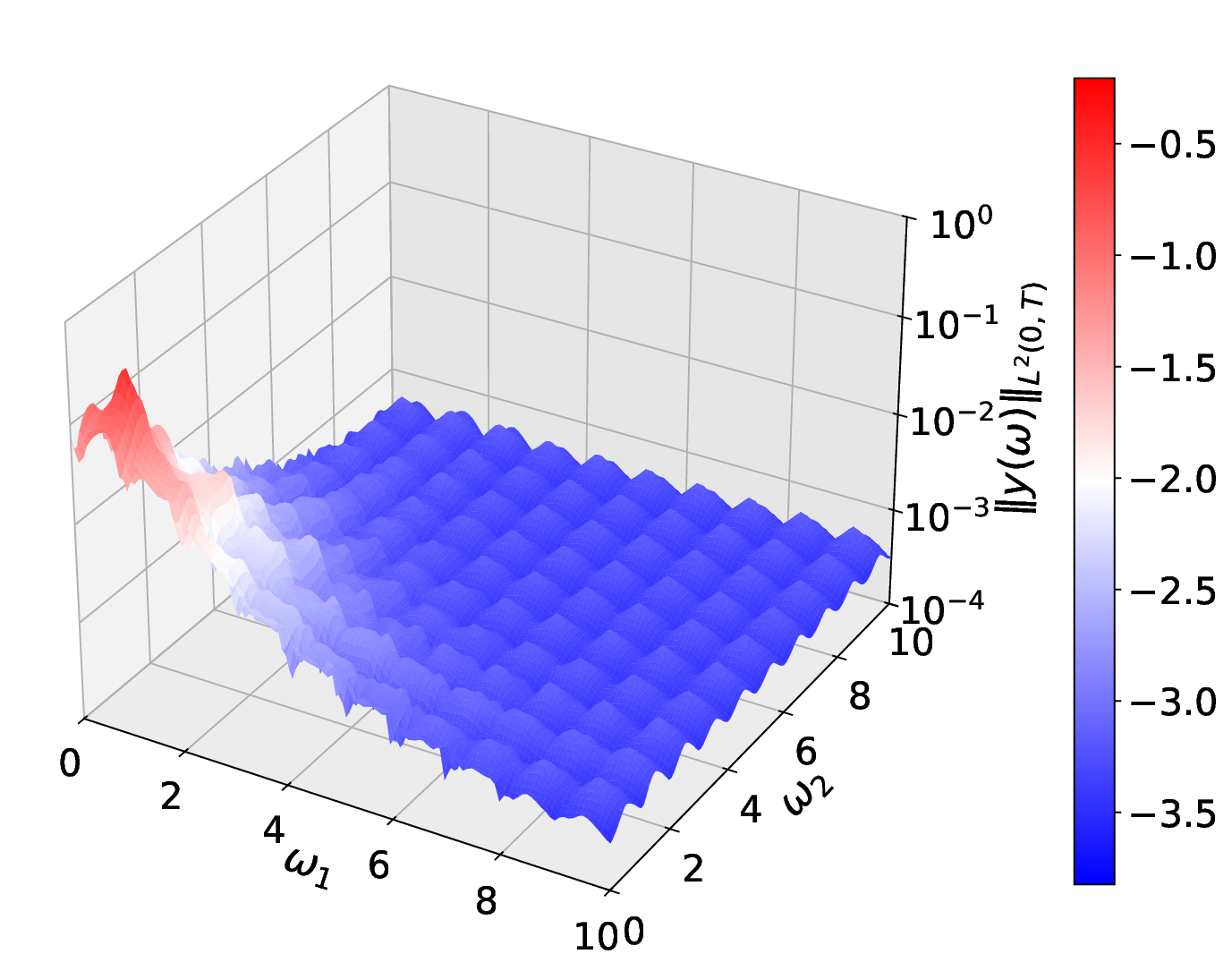}
    \caption{Spatial distribution of the temporal $L^2$-norm of the displacement for pure boundary (left) and grid control (right).}
    \label{Fig: L2Time}
\end{figure}

\begin{figure}[htbp]
    \centering
    \includegraphics[width=0.48\textwidth]{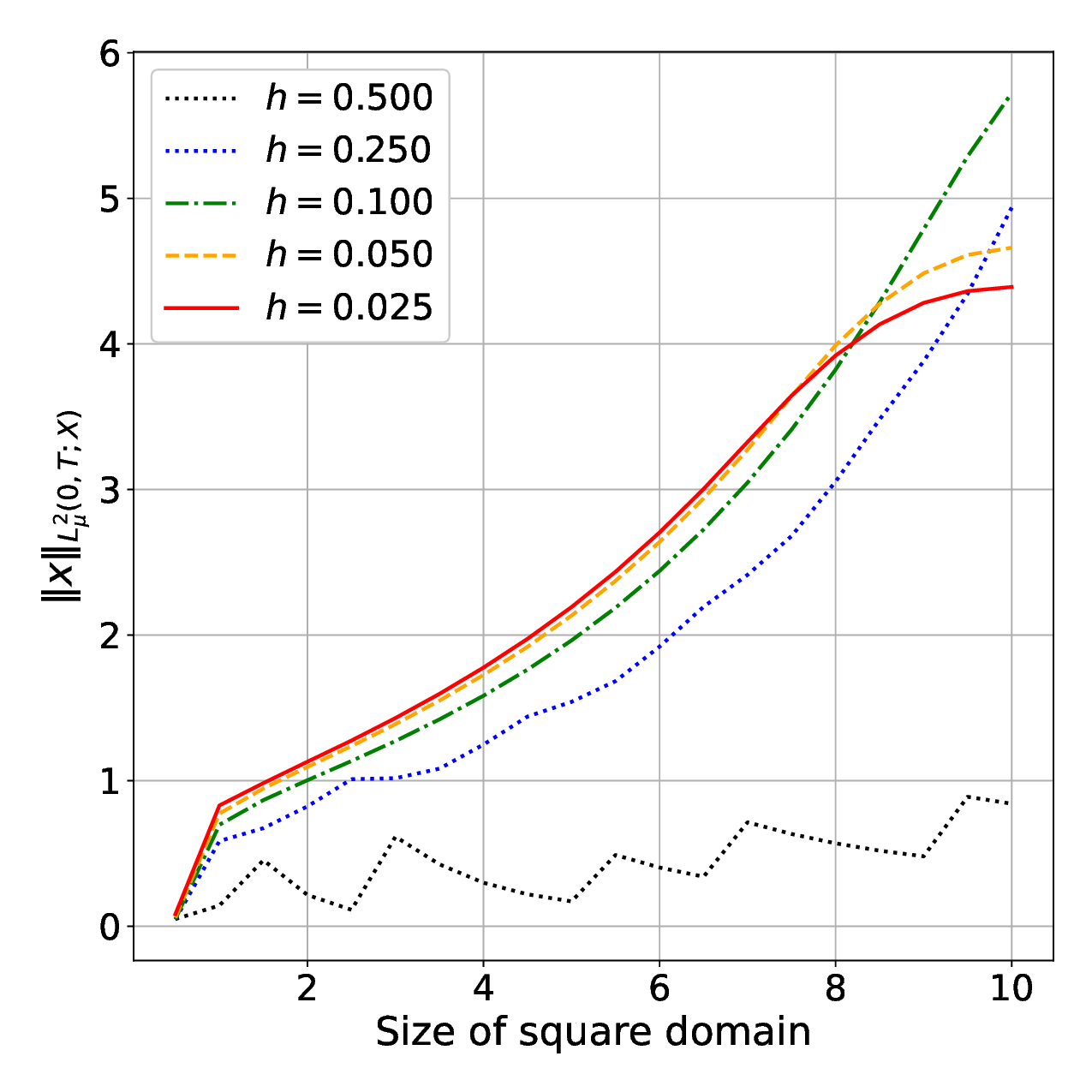}
    \hfill
    \includegraphics[width=0.48\textwidth]{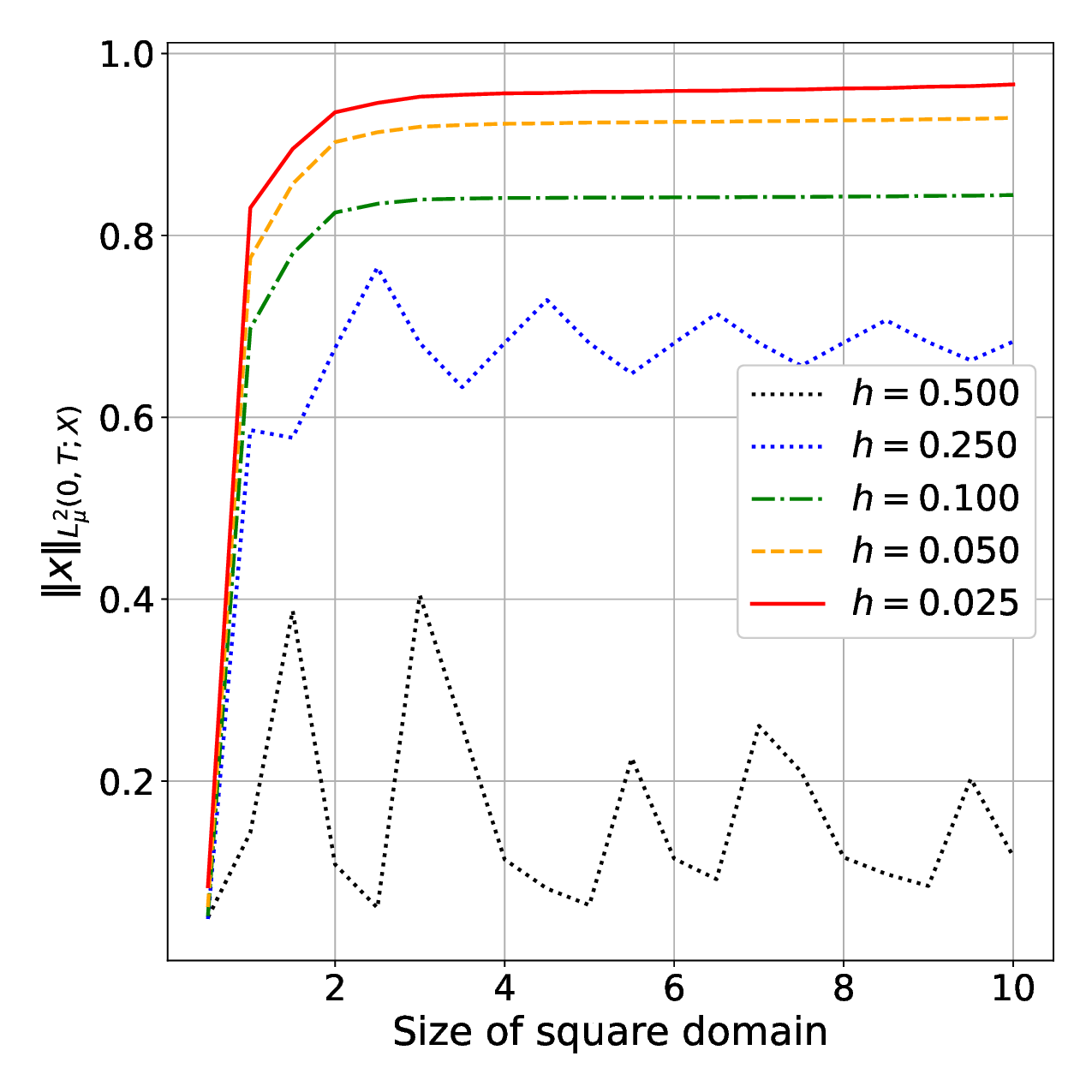}
    \caption{Comparison of an exponentially weighted norm for pure boundary (left) and grid control (right) for several grid resolutions.}
    \label{Fig: L2Weighted}
\end{figure}

\noindent The remaining question to be answered is, whether the decay under the influence of grid control which we noticed in Figures~\ref{Fig: Displacement}-\ref{Fig: Energy} is connected to exponential localization, or if we are just observing a higher order polynomial degree.
The following two Figures~\ref{Fig: L2Time} and~\ref{Fig: L2Weighted} provide strong evidence with regard to this subject.
\noindent \Cref{Fig: L2Time} shows a semilogarithmic plot of the distribution of the temporal $L^2$-norm
\begin{equation*}
    \n{y_\mathrm{opt}(\omega,\cdot)} _{L^2(0,T)}:= \int _0^T \n{y_\mathrm{opt}(\omega,t)}_2^2 \mathrm{d}t
\end{equation*}
over the spatial domain $\Omega = [0,10]\times [0,10]$. This norm is well defined since the displacement $y_\mathrm{opt}$ is weakly differentiable, i.e. its trace is well defined.
For Scenario 2, we observe a clear spatial decay as the distance from the initial perturbation increases. 
In the semilogarithmic representation, this decay is approximately linear over a substantial part of the domain, which indicates an exponential decay behaviour.

\noindent The exponential decay is, however, no longer resolved accurately once the values reach approximately $5\cdot 10^{-4}$. Below this level, the profile enters an uneven plateau rather than continuing to decay smoothly.
This effect is most likely caused by the finite spatial discretization, rather than being a feature of the optimally controlled wave equation. On the other hand no decay is visible on the semilogarithmic scale when examining the illustration of the norm for pure boundary control.
This observation is only consistent with our previous considerations, where noticed, that the wave is propagated across the whole spatial domain without being damped.
\noindent \Cref{Fig: L2Weighted} shows a comparison of the optimal boundary and grid control w.r.t. the exponentially weighted state space norm 
\begin{equation*}
    \n{x_\mathrm{opt}}_{X_\mu}
    :=\n{e^{\mu\n{P-\cdot}_2}y_\mathrm{opt}} _{L^2(0,T;H^1(\Omega))}
    + \n{e^{\mu\n{P-\cdot}_2}v_\mathrm{opt}} _{L^2(0,T;L^2(\Omega))}.
\end{equation*}
For this simulation the decay parameter was set to $\mu = 0.2$.
For scenario 1 we observe a steep nonlinear growth of the norm with the domain size for all resolutions except the very coarse grid size $h=0.4$.
On the other hand for scenario 2 the norm seems asymptotically approach an upper bound.
This boundedness from above is precisely what is required by exponential localisation.
This seems to support our hypothesis, that the impact of the perturbation is localized for scenario 1 but not for scenario 2.

\section{Conclusion}

\noindent In this paper, we investigated spatial localization effects for an optimally controlled two-dimensional wave equation.
In particular we tested whether a uniformly distributed lower-dimensional control geometry can induce localization behaviour similar to that predicted by the results in~\cite{Oppeneiger2026}.

\noindent The numerical results support this hypothesis. 
For pure boundary control, the perturbation is able to propagate over large parts of the spatial domain before interacting with the controlled boundary. 
As a consequence, the observed decay deteriorates as the domain size increases. 
In contrast, when the control acts on a rectangular grid of interior line interfaces (see~\Cref{fig: SpatialDomain}), the perturbation is damped much earlier, and the spatial profiles of the time-integrated displacement exhibit an exponential decay in the distance to the source of the perturbation. 

\noindent The computations were performed for several spatial resolutions in order to distinguish genuine localization effects from discretization artefacts. Nonetheless, the simulations should be interpreted as numerical evidence rather than as a proof of domain-uniform localization.

\noindent Overall, the results suggest that localization phenomena known for one-dimensional wave equations with distributed control may also appear in higher-dimensional wave equations with lower-dimensional control sets. 
Building on this important insight future research may successfully extend the results on domain-uniform stabilizability and detectability of hyperbolic equations to boundary and line controls in multiple space dimensions.

\section*{Acknowledgements}
\noindent The author thanks Manuel Schaller and Karl Worthmann for helpful comments and discussions.

\printbibliography

@article{Ascher1999,
	author = {{Ascher, Uri M.} and {Reich, Sebastian}},
	title = {The midpoint scheme and variants for Hamiltonian systems: advantages and pitfalls},
	journal = {SIAM Journal on Scientific Computing},
	year = 1999,
	volume = 21,
	number = 3,
	pages = "1045--1065",
}

@article{Bardos1992,
	author = {{Bardos, Claude} and {Lebeau, Gilles} and {Rauch, Jeffrey}},
	title = {Sharp sufficient conditions for the observation, control, and stabilization of waves from the boundary},
	journal = {SIAM Journal on Control and Optimization},
	year = 1992,
	volume = 30,
	number = 5,
	pages = "1024--1065",
}

@article{Benzi2015,
	author = {{Benzi, Michele} and {Simoncini, Valeria}},
	title = {Decay bounds for functions of Hermitian matrices with banded or Kronecker structure},
	journal = {SIAM Journal on Matrix Analysis and Applications},
	year = 2015,
	volume = 36,
	number = 3,
	pages = "1263--1282",
}

@article{Breiten2020,
	author = {{Breiten, Tobias} and {Pfeiffer, Laurent}},
	title = {On the turnpike property and the receding-horizon method for linear-quadratic optimal control problems},
	journal = {SIAM Journal on Control and Optimization},
	year = 2020,
	volume = 58,
	number = 2,
	pages = "1077--1102",
}

@article{Damm2014,
	author = {{Damm, Tobias} and {Gr{\"u}ne, Lars} and {Stieler, Marleen} and {Worthmann, Karl}},
	title = {An exponential turnpike theorem for dissipative discrete time optimal control problems},
	journal = {SIAM Journal on Control and Optimization},
	year = 2014,
	volume = 52,
	number = 3,
	pages = "1935--1957",
}

@article{Demko1984,
	author = {{Demko, Stanley} and {Moss, William F.} and {Smith, Philip W.}},
	title = {Decay rates for inverses of band matrices},
	journal = {Mathematics of Computation},
	year = 1984,
	volume = 43,
	number = 168,
	pages = "491--499",
}

@article{Faulwasser2022,
	author = {{Faulwasser, Timm} and {Gr{\"u}ne, Lars}},
	title = {Turnpike properties in optimal control: An overview of discrete-time and continuous-time results},
	journal = {Handbook of Numerical Analysis},
	year = 2022,
	volume = 23,
	pages = "367--400",
}

@article{Goettlich2025,
	author = {{G{\"o}ttlich, Simone} and {Schaller, Manuel} and {Worthmann, Karl}},
	title = {Perturbations in {PDE}-constrained optimal control decay exponentially in space},
	journal = {ESAIM: COCV},
	year = 2025,
	volume = 31,
	pages = "27",
}

@article{Oppeneiger2026,
	author = {{G{\"o}ttlich, Simone} and {Oppeneiger, Benedikt} and {Schaller, Manuel} and {Worthmann, Karl}},
	title = {Spatial exponential decay of perturbations in optimal control of general evolution equations},
	journal = {ESAIM: COCV},
	year = 2026,
	volume = 32,
	pages = "11",
}

@article{Schaller2020,
	author = {{Gr{\"u}ne, Lars} and {Schaller, Manuel} and {Schiela, Anton}},
	title = {Exponential sensitivity and turnpike analysis for linear quadratic optimal control of general evolution equations},
	journal = {Journal of Differential Equations},
	year = 2020,
	volume = 268,
	number = 12,
	pages = "7311--7341",
}

@article{Schaller2022,
	author = {{Gr{\"u}ne, Lars} and {Schaller, Manuel} and {Schiela, Anton}},
	title = {Efficient model predictive control for parabolic {PDE}s with goal oriented error estimation},
	journal = {SIAM Journal on Scientific Computing},
	year = 2022,
	volume = 44,
	number = 1,
	pages = "A471--A500",
}

@book{Hairer2006,
	author = {{Hairer, Ernst} and {Lubich, Christian} and {Wanner, Gerhard}},
	title = {Geometric Numerical Integration: Structure-Preserving Algorithms for Ordinary Differential Equations},
	publisher = {Springer},
	edition = 2,
	year = 2006,
}

@article{Shin2022,
	author = {{Shin, Sungho} and {Anitescu, Mihai} and {Zavala, Victor M.}},
	title = {Exponential decay of sensitivity in graph-structured nonlinear programs},
	journal = {SIAM Journal on Optimization},
	year = 2022,
	volume = 32,
	number = 2,
	pages = "1156--1183",
}

@article{Shin2023,
	author = {{Shin, Sungho} and {Zavala, Victor M.}},
	title = {Diffusing-horizon model predictive control},
	journal = {IEEE Transactions on Automatic Control},
	year = 2023,
	volume = 68,
	number = 1,
	pages = "188--201",
}
\end{document}